\documentclass[]{article}
\usepackage[letterpaper,margin=01.2 in]{geometry}
\usepackage[utf8]{inputenc}
\usepackage{amssymb}
\usepackage{amssymb, amsmath}
\usepackage{graphicx}
\usepackage{caption}
\usepackage{subcaption}
\newtheorem{theorem}{Theorem}
\newtheorem{lemma}[theorem]{Lemma}
\title{Haar wavelet method for the coupled degenerate reaction diffusion PDEs and the ODEs having non-linear source }
\author{Meena Pargaei and B.V. Rathish Kumar \\ Department of Mathematics and Statistics, Indian Institute of Technology, Kanpur}
\date{}
\begin{document}

\maketitle

\begin{abstract}
In this work we propose the Haar wavelet method for the coupled degenerate reaction diffusion PDEs and the ODEs having non-linear
source with Neumann boundary, applicable in various fields of the natural sciences,engineering 
and economics, for example in gas dynamics, certain biological models, assets pricing in economics, composite media etc. Convergence
analysis of the proposed numerical scheme has been carried out. We use the GMRES solver to solve the linear system of equations.
Numerical solutions for the model problems of medical significance have been successfully solved.
\end{abstract}

\section{Introduction}
Degenerate reaction diffusion system arises in the mathematical modeling of the various fields of the natural sciences,engineering 
and economics, for example in gas dynamics, certain biological models, assets pricing in economics, composite media etc. The 
degeneracy into the model is corresponding to the interface between the two separate medium in the physical problem. This types of 
problems are not only important from the application point of view but equally interesting from the analysis point, since it asks for
the design of techniques for the existence, uniqueness and stability of the solutions. In sedimentation processes and traffic flow 
problems, the concentration of the local solids is modeled by a strongly degenerate parabolic equation \cite{rd1}. 

Mathematically modeled reaction-diffusion equations describes the variation in the concentration of one or more substances in the
separated spaces with the influence of the local chemical reactions and the diffusion. This description implies that this types of
systems are applied in chemistry, however, this system can also describes the dynamical processes of the biology, geology, physics 
and finance. Mathematically, reaction-diffusion systems take the form of semi-linear parabolic partial differential equations
\cite{rd2}. The system corresponding to the population dynamics of the spruce band-worm for a non-degenerate case in a biological 
setting is discussed in 
\cite{rd3}. On the other hand, similar governing equations also arises in mathematical biology as a well-known reaction-diffusion 
system modeling the interaction between two chemical species. Under certain conditions, it produces stationary solutions with Turing-type spatial patterns \cite{rd3,rd5} and a standard proof for the existence and uniqueness can be found in \cite{rd4}.
The difficulty in direct measurement of the cardiac electric activity motivates for the mathematical modeling and numerical 
simulations of this phenomena. Hodking and Huxely in 1952 modeled the first mathematical model to calculate the action potential 
in a squid giant axon which later modified to describe the several biological phenomena. Tung \cite{rd6} introduced the first 
mathematical model, called as Bidomain model
, for the study of the cardiac electric activity. This model consists of the two degenerate 
parabolic reaction diffusion system corresponding to the two spaces separated by the
interface membrane. This degenerate structure of the bidomain model is essentially due to the differences between the intra- and extracellular anisotropy of the cardiac tissue. 
 Colli Franzone and Savar \cite{rd7} present a weak formulation for the bidomain model and show that it has a structure suitable to
apply the theory of evolution
variational inequalities in Hilbert spaces. Bendahmane and Karlsen \cite{rd9} prove existence and uniqueness for
the bidomain equations using, for the existence part, the Faedo-Galerkin method and compactness theory,
and Bourgault, Coudi`ere, and Pierre \cite{rd9} prove existence and uniqueness for the bidomain equations, first 
reformulating the problem into a single parabolic PDE and then applying a semigroup approach. In [] Galerkin finite 
element error analysis for the coupled nonlinear degenerate system of advection - diffusion equations modeling a 
two-phase immiscible flow through porous media is derived. 

From a computational point of view, this space - time bidomain model has been numerically solved via finite 
difference method, finite volume method, finite element method, adaptive finite element methods using a posteriori 
error techniques, domain decomposition method using an alternating direction implicit method in []. In [], 
multiresolution technique is used to solve the degenerate system including the monodomain and bidomain models. 

Recently wavelets are getting much attention for their effective use to numerically solve the mathematical models from various field of science, engineering and biology. 
Wavelets are well known for their inherent nature to adopt to the complexities such as discontinuities, sharp variation etc.
Its properties such as orthogonality, compact support, arbitrary regularity and high order
vanishing moments are very attractive. Because of these properties, solutions with discontinuities or fast oscillations 
in a localized region, can be approximated well by using very few wavelets. This method has been used to find the solution of the Integral equations, ordinary differential equations, partial differential equations,
and fractional partial differential equations \cite{aziz,siraj,siraj2,aziz2,siraj3,siraj4,aziz3,
siraj5,naldi,wu}. Different kinds of popular wavelets such as Daubechies \cite{diaz}, Chebyshev \cite{Babolian}, 
Haar wavelets \cite{lpeik},
Battle-Lemarie \cite{zhu}, B-spline \cite{Dehghan}, Legendre wavelets \cite{wang} are being used by the researchers for different 
models. Out of these wavelets Haar wavelets are very popular because of its simplicity and easy implementation in finite domains. 
Haar wavelets are piecewise constant functions which are orthogonal and have compact support. They also have scaling property.
Because of these properties of haar wavelets, Haar wavelet method has become very popular. In this paper, we will discuss the Haar 
wavelets and the collocation based haar wavelet method.
Because of the discontinuity of Haar wavelets, the derivatives does not exist. In this situation it is not possible to calculate
the solution of differential equations. To overcome this difficulty,  Chen and Hsiao \cite{chen} have
proposed an idea that the highest order derivative of the differential equation is expanded into 
the Haar series, not the function itself. Then on integration one can obtain lower order derivatives and the 
functions too. 

Haar wavelet method has been used to solve the linear and non-linear ordinary differential equations of
all order in \cite{bujurke}, elliptic and parabolic partial differential equations with Dirichlet and Neumann 
boundary conditions both in \cite{lpeik, kannan, sheo} and also for the eigenvalue problems in \cite{ev}.  The basic
technique of the collocation Haar wavelet is to convert the continuous problem into a discrete form with finite number 
of collocation points. Haar wavelet has been used vastly in the field of signal processing 
communication, Image processing. In \cite{siraj2}, fluid flow boundary layer problem is solved via Haar wavelet
collocation method.  Haar wavelet method has also been used to solve the nonlocal problem in 
two-dimension. I. Singh and S. Kumar \cite{sheo} proposed the Haar wavelet collocation method for the solution of three 
dimensional Poisson and Helmholtz equations. Haar wavelet method has been used to solve the wave-like equations by 
B. Naresh et.al \cite{naresh}.
In this work, we will develop the haar wavelet method for the coupled non-linear degenerate PDEs-ODEs system with Neumann boundary condition. We will show the advantage of haar wavelet method, like easy implementation and easily 
extendable to the higher dimension.
 In the next section, we will introduce the haar wavelet function, properties, and 
their integration functions. One-dimensional, two dimensional and the three dimensional haar wavelet method for the
non-linear coupled non-linear degenerate PDE-PDE system with Neumann boundary which is coupled with the system of ODEs will be developed in section 3. Convergence analysis of the proposed method is 
conducted in section 4. Numerical results and discussion for problems in one, two and three dimensions has been discussed in the next section.\\

\section{Haar wavelets}
Let us consider the interval $x \in [A,B]$, A, B are finite real numbers. Define $M=2^J$, $J$ is 
the maximum level of resolution. This interval $[A,B]$ is equally divided into $2M$ subintervals 
such that the length of each subinterval is $\Delta x = (B-A)/2M$. Now, define the dilation and 
translation parameter $j=0,1,...,J$ and $k=0,1,...,m-1$ respectively, where, $m=2^j$. The wavelet number is given by $i=m+k+1$. 
Family of haar wavelets is defined as follows:

For $i \geq 2$
\begin{align}
\label{haar}
h_i(x) = 
\begin{cases}
   1 & \text{when } \beta_1(i) \leq \ x < \beta_2(i)\\    
 - 1 & \text{when } \beta_2(i) \leq \ x < \beta_3(i)\\  
   0 & \text{otherwise},
\end{cases}
\end{align}

where 

\begin{align*}
 & \beta_1(i)= A+ 2k\zeta \Delta x, \hspace{0.5cm} \beta_2(i)= A+ (2k+1)\zeta \Delta x, \\
 & \beta_3(i)= A+ 2(k+1)\zeta \Delta x,  \hspace{0.5cm} \zeta=M/m.
\end{align*}

For $i=1$
\begin{align*}
h_i(x) = 
\begin{cases}
   1 & \text{when } A \leq \ x < B\\    
   0 & \text{otherwise}.
\end{cases}
\end{align*}

Haar wavelets are orthogonal, since
\begin{align*}
\int_A^B h_i(x) h_j(x) dx = 
\begin{cases}
   2^{-j}(B-A) & \text{when } i=j\\    
   0 & \text{when } i\neq j.
\end{cases}
\end{align*}

For the solution of differential equation, we have to compute the integral
\begin{align}
 p_{\alpha,i}(x)= \idotsint_A^B h_i(t)dt^{\alpha} = \frac{1}{(\alpha -1)!}
 \int_A^x (x-t)^{(\alpha -1)} h_i(t) dt,
\end{align}

where $\alpha = 1, 2,..., n$ and $i=1,2,..., 2M$. 

For the case $\alpha =0, p_{0,i}(x)=h_i(x)$

This integral is calculated with the help of equation \eqref{haar}, which is given by

\begin{align}
\label{HaarIntegral}
p_{\alpha ,i}(x) = 
\begin{cases}
   0 & \text{when } x < \beta_1(i)\\    
 \frac{1}{\alpha !}(x-\beta_1(i))^{\alpha} & \text{when } \beta_1(i) \leq \ x < \beta_2(i)\\ 
  \frac{1}{\alpha !}[(x-\beta_1(i))^{\alpha}-2(x-\beta_2(i))^{\alpha}] & \text{when } \beta_2(i) \leq \ x < \beta_2(i)\\  
   \frac{1}{\alpha !}[(x-\beta_1(i))^{\alpha}-2(x-\beta_2(i))^{\alpha}+(x-\beta_3(i))^{\alpha})] & \text{when }  x > \beta_3(i). 
\end{cases}
\end{align}

When $\alpha =1, 2$,
\begin{align}
p_{1,i}(x) = 
\begin{cases}
   (B-A) & \text{when } i=1\\    
   0 & \text{otherwise}.
\end{cases}
\end{align}

\begin{align*}
p_{2,i}(x) = 
\begin{cases}
   (B-A)^2/2 & \text{when } i=1\\    
   (B-A)^2/4m^2 & \text{otherwise}.
\end{cases}
\end{align*}

For the grid points $x_k=A+kh, k=0,1,...,2M$, $h=\Delta x$, collocation points
are as follows:
\begin{align}
 y_k= \frac{x_{k-1}+x_{k}}{2}, k=1,2,...,2M.
\end{align}

After this discretization, we define Haar matrix $H$, and Haar Integral matrices $P_1 , P_2$ of
size $2M \times 2M$ as $H(i,k)=h_i(y_k), P_1(i,k)=p_{1,i}(y_k) , P_2(i,k)=p_{2,i}(y_k)$.

\subsection{Function approximation}
Any function $u(x) \in L^2[0,1)$ can be approximated in terms of the Haar wavelet series as 

\begin{align*}
f(x)=\sum_{i=1}^{\infty}\alpha_i h_i(x),
\end{align*}

where the wavelet coefficients $\alpha_i$ are obtained by

\begin{align*}
\alpha_i= 2^j \int_0^1 f(x) h_i(x)dx.
\end{align*}

Since only the finite number of terms are taken for the computational purpose therefore, 
the function approximation $f$ is given by
 \begin{align*}
 f(x)=\sum_{i=1}^{n}\alpha_i h_i(x).
 \end{align*}

\section{Mathematical Model}
Consider the degenerate parabolic reaction diffusion system coupled with a system of ODEs of the type
\begin{align}
\nonumber
& \frac{\partial v}{\partial t}- div(D_i(x)\nabla u_i) + f(v,w)=I_1, \hspace{1cm} x \in \Omega , t \in (0,T) \\ \nonumber
&\frac{\partial v}{\partial t}- div(D_e(x)\nabla u_e) + f(v,w)= I_1, \hspace{1cm} x \in \Omega, t \in (0,T)\\ \nonumber
&\frac{\partial w}{\partial t}-g(v,w)=0, \hspace{1cm} x \in  \Omega,  t \in (0,T)\\ \nonumber
&v(x,0)= v_0(x,0), \hspace{5mm} w(x,0)=w_0(x,0) \hspace{1.8cm}  x \in \Omega\\ \nonumber
&n^T D_{i,e}(x) \nabla v =0 \hspace{5.1cm}  x \in \partial \Omega, t \in (0,T). 
\end{align}
where, $v=u_i-u_e$. 
We can rewrite the above set of equations as follows:
\begin{align}
 \nonumber
&C_m \frac{\partial v}{\partial t}+div(D_e(x)\nabla u_e) + f(v,w)= I_2, \hspace{1.5cm} x \in \Omega , t \in (0,T) \\  \nonumber
&- div((D_i+)D_e)(x)\nabla u_e)-div(D_i(x)\nabla v)=I_1-I_2, \hspace{0.6cm} x \in \Omega, t \in (0,T)\\  \nonumber
&\frac{\partial w}{\partial t}-g(v,w)=0, \hspace{1cm} x \in  \Omega,  t \in (0,T)\\  \nonumber
&v(x,0)= v_0(x,0), \hspace{5mm} w(x,0)=w_0(x,0) \hspace{1.8cm}  x \in \Omega\\  \nonumber
&n^T D_{i,e}(x) \nabla v =0 \hspace{ 5.2cm}  x \in \partial \Omega, t \in (0,T). 
\end{align}

\subsection{Haar wavelet method for the coupled degenerate reaction diffusion PDE and the ODEs}

Consider the following system of coupled degenerate PDEs and ODEs in one dimension:
\begin{align}
\label{mvb1d}
& \frac{\partial v}{\partial t}-\frac{d}{dx}(D_e(x)\frac{d}{dx} u_e) + f(v,w)=I_2, \hspace{1.5cm} x \in \Omega , t \in (0,T) \\ 
\label{mueb1d}
&- \frac{d}{dx}((D_i+)D_e)(x)\frac{d}{dx} u_e)-\frac{d}{dx}(D_i(x)\frac{d}{dx} v)=I_1-I_2, \hspace{0.6cm} x \in \Omega, t \in (0,T)\\
\label{mwb1d}
&\frac{\partial w}{\partial t}-g(v,w)=0, \hspace{1cm} x \in  \Omega,  t \in (0,T)\\
&v(x,0)= v_0(x,0), \hspace{5mm} w(x,0)=w_0(x,0), \hspace{1.8cm}  x \in \Omega\\
&D(x)\frac{d u_{i,e}}{d x}(0,t)=0, \hspace{0.5cm} D(x)\frac{du_{i,e}}{d x}(1,t)=0 \hspace{1cm}  x \in \partial \Omega, t \in (0,T). 
\end{align}

Let us take
\begin{align}
\label{discvb}
&\frac{\partial^3 v}{\partial t \partial x^2}(x,t) = \sum_{i=1}^{2*M}\alpha_i h_i(x), \hspace{1cm}t\in [t_s, t_{s+1}),\\
\label{discueb}
&\frac{\partial^2 u_e}{\partial x^2}(x,t) = \sum_{j=1}^{2*M}\beta_j h_j(x), \hspace{1.2cm}t\in [t_s, t_{s+1}),\\
\label{discwb}
&\frac{\partial w}{\partial t}(x,t)= \sum_{m=1}^{2*M}\gamma_m h_m(x), \hspace{1.2cm}t\in [t_s, t_{s+1}).
\end{align}
Integrating equation \eqref{discvb} w.r.t. $t$ from $t_s$ to $t$, we obtain
\begin{align}
\label{d2vb}
& \frac{\partial^2 v}{\partial x^2}(x,t)= (t-t_s) \sum_{i=1}^{2M}\alpha_i h_i(x)+ \frac{\partial^2 v}{\partial x^2}(x,t_s).
\end{align}
Integrating equation \eqref{d2vb} w.r.t. $x$ from $0$ to $x$ twice and using the boundary condition on $v$, we obtain

\begin{align}
\label{dv1b}
& \frac{\partial v}{\partial x}(x,t)= (t-t_s) \sum_{i=1}^{2M}\alpha_i P_{1,i}(x)+ \frac{\partial v}{\partial x}(x,t_s)\\
\label{vb}
& v(x,t)=(t-t_s) \sum_{i=1}^{2M}\alpha_i P_{2,i}(x)+ v(x,t_s)-v(0,t_s)+v(0,t).
\end{align}
Again, Integrating equation \eqref{discueb} w.r.t. $x$ from $0$ to $x$ twice and using the boundary condition on $u_e$, we obtain
\begin{align}
\label{due1b}
& \frac{\partial u_e}{\partial x}(x,t)= \sum_{j=1}^{2M}\beta_j P_{1,j}(x)+ \frac{\partial u_e}{\partial x}(0,t_s),\\
\label{ueb}
& u_e(x,t)=(t-t_s) \sum_{j=1}^{2M}\beta_j P_{2,j}(x)+ u_e(0,t).
\end{align}

Now, Integrating equation \eqref{discvb} w.r.t. $x$ from $0$ to $x$ twice and using the boundary condition on $v$, we obtain

\begin{align}
\label{dvtb}
& \frac{\partial^2 v}{\partial t \partial x}(x,t)= (t-t_s) \sum_{i=1}^{2M}\alpha_i P_{1,i}(x)+ \frac{\partial^2 v}{\partial t \partial x}(0,t),\\
& \frac{\partial v}{\partial t}(x,t)= (t-t_s) \sum_{i=1}^{2M}\alpha_i P_{2,i}(x)+ \frac{\partial v}{\partial x}(0,t).
\end{align}
 To calculate the solution of the system \eqref{mvb1d}-\eqref{mwb1d} at the grid points we will write it in the discrete form as 
 follows: 
\begin{align}
\label{discmvb1d}
& \frac{\partial v}{\partial t}(x_k,t_{s+1})-D_e(x)\frac{d^2 v}{d x^2})(x_k,t_{s+1})- D_e'(x)\frac{d v}{d x})(x_k,t_{s+1})
+f(v(x_k,t_{s}),\\ \nonumber
&w(x_k,t_{s+1})) = 0,\\
& -(D_i(x)+D_e(x))\frac{d^2 v}{d x^2})(x_k,t_{s+1})- \big (D_i'(x)+D_e'(x)\big )\frac{d v}{d x})(x_k,t_{s+1}) = -I,\\
\label{discmwb1d}
& \frac{\partial w}{\partial t}(x_k,t_{s+1})=g(v(x_k,t_{s+1}),w(x_k,t_{s+1}).
\end{align}
Now, using \eqref{d2vb},\eqref{dv1b} in \eqref{discmvb1d}, we will get,
\begin{align}
\label{mat1}
\nonumber
&\sum_{i=1}^{2M}\alpha_i P_{2,i}(x_k) +  \sum_{j=1}^{2M} \beta_j \big [D_e(x_k) h_j(x_k)+D_e'(x_k)P_{1,j}(x_k)\big ] = 
D_e'(x_k)u_e(x_1,t_{s+1})\\ 
& +\frac{\partial v}{\partial t}(x_1,t_{s+1}) -f(v(x_k,t_s), w(x_k,t_{s+1})).
\end{align}
Similarly, using \eqref{mueb1d}, \eqref{due1b} in \eqref{discmwb1d}, we obtain,
\begin{align}
\label{mat2}
\nonumber
&\sum_{i=1}^{2M}\alpha_i \big [D_i(x_k)(t_{s+1}-t_s)h_i(x_k) + D_i'(x_k)(t_{s+1}-t_s)P_{1,i}(x_k)] + \sum_{j=1}^{2M}\beta_j \big [(D_i(x_k)+ \\ \nonumber
 & D_e(x_k)) h_j(x)+(D_i(x_k)+ D_e(x_k))' P_{1,j}(x_k) \big ] =
D_i(x_k)\frac{d^2 v}{d x^2}(x_k,t_{s})- D_i'(x_k)\frac{d v}{d x}\\ &(x_k,t_{s}) + I.
\end{align}
Now, from \eqref{mat1} and \eqref{mat2}, we will get the following matrix system,
\begin{align}
\label{matv}
K \begin{bmatrix} \alpha & \beta \end{bmatrix} = b,
\end{align}

Now, Integrate \eqref{discwb} w.r.t. $t$ from $t_s$ to $t$, we get
\begin{align}
w(x,t)=(t-t_s) \sum_{i=1}^{2M}\gamma_m h_m(x)+ w(x,t_s)
\end{align}
Using \eqref{discwb} in equation \eqref{discmwb1d}, and linearize the non-linear term taking values at the previous time step, we get
\begin{align*}
\sum_{m=1}^{2M}\gamma_m h_m(x_k) &= g(v(x_k,t_{s}),w(x_k,t_{s})).
\end{align*}

Matrix system of the above equation is given by,
\begin{align}
\label{matw}
H\gamma=c,
\end{align}
At each time step, firstly we will calculate the $w$ at the desired time by solving equation \eqref{matw} and then obtain $v$ using
obtained $w$, at the desired time.

\subsection{Haar wavelet method for the coupled degenerate reaction diffusion PDE and the ODEs in two dimension}

Consider the coupled degenerate reaction diffusion PDEs and the ODEs given as follows:
\begin{align}
\label{mvb2d}
C_m \frac{\partial v}{\partial t}- \nabla. (D_e(x,y)\nabla u_e) + f(v,w) & =I_1, \hspace{1.5cm}  0\leq x,y \leq 1,  0\leq t \leq T\\ 
\label{mueb2d}
- \nabla ((D_i+D_e)(x,y)\nabla u_e)-\nabla (D_i(x,y)\nabla v)&=I_1-I_2, \hspace{0.6cm} 0\leq x,y \leq 1,  0\leq t \leq T\\
\label{mwb2d}
\frac{\partial w}{\partial t}-g(v,w) &=0, \hspace{1cm} 0\leq x,y \leq 1,  0\leq t \leq T\\
v(x,y,0)= v_0(x,y,0), \hspace{5mm} w(x,y,0) &=w_0(x,0), \hspace{1.8cm} 0\leq x,y \leq 1\\
D(x,y)\frac{d u_{i,e}}{d x}(0,y,t)&=0,\hspace{2cm} 0\leq t \leq T\\
D(x,y)\frac{d u_{i,e}}{d x}(1,y,t) &=0 \hspace{2cm}, 0\leq t \leq T\\
D(x,y)\frac{d u_{i,e}}{d y}(x,0,t)&=0,\hspace{2cm} 0\leq t \leq T\\
D(x,y)\frac{d u_{i,e}}{d y}(x,1,t) &=0, \hspace{2cm}, 0\leq t \leq T
\end{align}

Let us write $\frac{\partial ^5 v}{\partial t \partial x^2 \partial y^2}(x,y,t)$, $\frac{\partial ^4 ue}{\partial x^2 \partial y^2}
(x,y,t)$
and $\frac{\partial w}{\partial t}(x,y,t)$ in terms of the Haar wavelet as follows :
\begin{align}
\label{discvb2d}
&\frac{\partial^5 v}{\partial t \partial x^2 \partial y^2}(x,y,tt) = \sum_{i,j=1}^{2*M}\alpha_{i,j} h_i(x)h_j(y), \hspace{1cm}t\in [t_s, t_{s+1}),\\
\label{discueb2d}
&\frac{\partial^4 u_e}{\partial x^2\partial y^2 }(x,t) = \sum_{m,n=1}^{2*M}\beta_{m,n} h_m(x)h_n(y), \hspace{1.2cm}t\in [t_s, t_{s+1}),\\
\label{discwb2d}
&\frac{\partial w}{\partial t}(x,t)= \sum_{r,s=1}^{2*M}\gamma_{r,s} h_r(x)h_s(y), \hspace{1.2cm}t\in [t_s, t_{s+1}).
\end{align}
Integrating equation \eqref{discvb2d} w.r.t  $t$ from $t_s$ to $t$, we will get 
\begin{align}
\label{vb4}
\frac{\partial ^4 v}{ \partial x^2 \partial y^2} (x,y,t) &= (t-t_s)\sum_{i,j=1}^{2M}\alpha_{i,j} h_i(x) h_j(y)
+\frac{\partial ^4 v}{ \partial x^2 \partial y^2}(x,y,t_s), \hspace{0.5cm}t\in [t_s, t_{s+1}).
\end{align}

Now, Integrate equation \eqref{vb4} twice w.r.t  $x$ from $0$ to $x$ also using the boundary conditions,
we will obtain the following
\begin{align}
&\frac{\partial ^3 v}{ \partial x \partial y^2}(x,y,t)= (t-t_s)\sum_{i,j=1}^{2M}\alpha_{i,j} P_{1,i}(x) h_j(y)
+\frac{\partial ^3 v}{ \partial x \partial y^2}(x,y,t_s), \hspace{0.5cm}t\in [t_s, t_{s+1}),\\ \nonumber 
\label{vb2y}
&\frac{\partial ^2 v}{\partial y^2}(x,y,t)= (t-t_s)\sum_{i,j=1}^{2M}\alpha_{i,j} P_{2,i}(x) h_j(y)
+\frac{\partial ^2 v}{\partial y^2}(x,y,t_s)- \frac{\partial ^2 v}{\partial y^2}(0,y,t_s)\\
&+\frac{\partial ^2 v}{\partial y^2}(0,y,t), \hspace{0.5cm}t\in [t_s, t_{s+1}).
\end{align}

Similarly, Integrate equation \eqref{vb4} twice w.r.t  $y$ from $0$ to $y$ also using the boundary conditions, we get
\begin{align}
&\frac{\partial ^3 v}{ \partial x^2 \partial y} (x,y,t)= (t-t_s)\sum_{i,j=1}^{2M}\alpha_{i,j} h_i(x) P_{1,j}(y) 
+\frac{\partial ^3 v}{ \partial x^2 \partial y}(x,y,t_s), \hspace{0.5cm}t\in [t_s, t_{s+1}),\\ \nonumber
\label{vb2x}
&\frac{\partial ^2 v}{\partial x^2}(x,y,t)= (t-t_s)\sum_{i,j=1}^{2M}\alpha_{i,j} h_i(x) P_{2,j}(y) 
+\frac{\partial ^2 v}{\partial x^2}(x,y,t_s)- \frac{\partial ^2 v}{\partial x^2}(x,0,t_s) \\
&+\frac{\partial ^2 v}{\partial x^2}(x,0,t), \hspace{0.5cm}t\in [t_s, t_{s+1}).
\end{align}

Now, Integrate \eqref{vb2x} w.r.t  $x$ from $0$ to $x$ and Integrate \eqref{vb2y} w.r.t  $y$ from $0$ to $y$ also 
using the boundary conditions,we get
\begin{align}
\label{vb1x}
\nonumber
&\frac{\partial v}{\partial x} (x,y,t)= (t-t_s)\sum_{i,j=1}^{2M}\alpha_{i,j} P_{1,i}(x) P_{2,j}(y) 
+\frac{\partial v}{\partial x}(x,y,t_s)- \frac{\partial v}{\partial x^2}(x,0,t_s)\\
&+\frac{\partial v}{\partial x}(x,0,t), \hspace{0.5cm}t\in [t_s, t_{s+1}).\\
\label{vb1y}
\nonumber
&\frac{\partial v}{\partial y}(x,y,t)= (t-t_s)\sum_{i,j=1}^{2M}\alpha_{i,j} P_{2,i}(x) P_{1,j}(y)
+\frac{\partial v}{\partial y}(x,y,t_s)-\frac{\partial v}{\partial y}(0,y,t_s) \\
&+\frac{\partial v}{\partial y}(0,y,t), \hspace{0.5cm}t\in [t_s, t_{s+1}).
\end{align}

Again, Integrating \eqref{vb1x} w.r.t  $x$ from $0$ to $x$, we obtain the following
\begin{align}
\nonumber
&v(x,y,t)=(t-t_s)\sum_{i,j=1}^{2M}\alpha_{i,j} P_{2,i}(x) P_{2,j}(y)+v(x,y,t_s)-v(0,y,t_s)-v(x,0,t_s)\\
&+v(0,0,t_s)+v(x,0,t)-v(0,0,t)+v(0,y,t).
\end{align}
Now, Integrate \eqref{discvb2d} twice w.r.t  $x$ from $0$ to $x$ and apply the boundary conditions, we get
\begin{align*}
&\frac{\partial ^3 v}{\partial y^2 \partial t} (x,y,t)= \sum_{i,j=1}^{2M}\alpha_{i,j}P_{2,i}(x) h_j(y) 
+\frac{\partial ^3 v}{\partial y^2 \partial t}(0,y,t).
\end{align*}

Now, Integrate above equation twice w.r.t $y$ from $0$ to $y$ and apply the boundary conditions, we obtain
\begin{align}
\nonumber
\label{vb1t}
&\frac{\partial v}{\partial t}(x,y,t)= \sum_{i,j=1}^{2M}\alpha_{i,j} P_{2,i}(x) P_{2,j}(y)
+\frac{\partial v}{\partial t}(0,y,t)- \frac{\partial v}{\partial t}(0,0,t)
+\frac{\partial v}{\partial t}(x,0,t),\\
&\hspace{9cm}t\in [t_s, t_{s+1}).
\end{align}

Now, Integrate equation \eqref{discueb2d} twice w.r.t  $x$ from $0$ to $x$ also using the boundary conditions,
we will obtain the following
\begin{align}
\nonumber
&\frac{\partial ^3 u_e}{ \partial x \partial y^2}(x,y,t)= \sum_{m,n=1}^{2M}\beta_{m,n} P_{1,m}(x) h_m(y),\\
\label{ueb2y}
&\frac{\partial ^2 u_e}{\partial y^2}(x,y,t)= \sum_{m,n=1}^{2M}\beta_{m,n} P_{2,m}(x) h_n(y)
+\frac{\partial ^2 u_e}{\partial y^2}(0,y,t)\hspace{0.5cm}t\in [t_s, t_{s+1}).
\end{align}

Similarly, Integrate equation \eqref{discueb2d} twice w.r.t  $y$ from $0$ to $y$ also using the boundary conditions, we get
\begin{align}
\nonumber
&\frac{\partial ^3 u_e}{ \partial x^2 \partial y} (x,y,t)=\sum_{m,n=1}^{2M}\beta_{m,n} h_m(x) P_{i,n}(y), \\
\label{ueb2x}
&\frac{\partial ^2 u_e}{\partial x^2}(x,y,t)= \sum_{m,n=1}^{2M}\beta_{m,n} h_m(x) P_{2,n}(y) 
+\frac{\partial ^2 u_e}{\partial x^2}(x,0,t) \hspace{0.5cm}t\in [t_s, t_{s+1}).
\end{align}
Now, Integrate \eqref{ueb2x} w.r.t  $x$ from $0$ to $x$ and Integrate \eqref{ueb2y} w.r.t  $y$ from $0$ to $y$ also 
using the boundary conditions,we get
\begin{align}
\label{ueb1x}
&\frac{\partial u_e}{\partial x} (x,y,t)= \sum_{m,n=1}^{2M}\beta_{m,n} P_{1,m}(x) P_{2,n}(y)+\frac{\partial u_e}{\partial x}(x,0,t), \hspace{0.5cm}t\in [t_s, t_{s+1}),\\
\label{ueb1y}
&\frac{\partial u_e}{\partial y}(x,y,t)= \sum_{m,n=1}^{2M}\beta_{m,n} P_{2,m}(x) P_{1,n}(y)
+\frac{\partial u_e}{\partial y}(0,y,t), \hspace{0.5cm}t\in [t_s, t_{s+1}).
\end{align}

Again, Integrating \eqref{ueb1x} w.r.t  $x$ from $0$ to $x$, we obtain the following
\begin{align}
\label{ue2d}
&u_e(x,y,t)=\sum_{m,n=1}^{2M}\beta_{m,n} P_{2,m}(x) P_{2,n}(y)+u_e(0,y,t)+u_e(x,0,t)-u_e(0,0,t).
\end{align}

Again, Integrate \eqref{discwb2d} w.r.t $t$ from $t_s$ to $t$, we acquire
\begin{align}
\label{wb2d}
w(x,y,t) = (t-t_s)\sum_{r,s=1}^{2M}\beta_{r,s} h_{r}(x) h_{s}(y)+ w(x,y,t_s)
\end{align}

Now, to calculate the solution of the system \eqref{mvb2d}-\eqref{mwb2d} at the grid points we will write it in the discrete 
form as follows: 
\begin{align}
\nonumber
\label{discmvb2d}
& \frac{\partial v}{\partial t}(x_k,y_l,t_{s+1})- \bigg[(\sigma _{e,l}(x_k,y_l)\frac{\partial ^2 ue}{\partial x^2}(x_k,y_l,t_{s+1})+ 
\sigma _{e,l,x}(x_k,y_l)\frac{\partial  ue}{\partial x}(x_k,y_l,t_{s+1})) + \\ \nonumber
&((\sigma _{e,t}(x_k,y_l)\frac{\partial ^2 ue}{\partial y^2}(x_k,y_l,t_{s+1})
 +\sigma _{e,t,y}(x_k,y_l)\frac{\partial ue}{\partial y}(x_k,y_l,t_{s+1})\bigg ]
+f(v(x_k,y_l,t_{s+1}), \\
& w(x_k,y_l,t_{s+1})) = I_2,\\
\nonumber
\label{discmueb2d}
&- \bigg[(\sigma _{i,l}+\sigma _{e,l})(x_k,y_l)\frac{\partial ^2 ue}{\partial x^2}(x_k,y_l,t_{s+1})+ 
(\sigma _{i,l,x}+\sigma _{e,l,x})(x_k,y_l)\frac{\partial  ue}{\partial x}(x_k,y_l,t_{s+1})) + \\ \nonumber
&((\sigma _{i,l}+\sigma _{e,l})(x_k,y_l)\frac{\partial ^2 ue}{\partial y^2}(x_k,y_l,t_{s+1})
 +(\sigma _{i,l,y}+\sigma _{e,l,y})(x_k,y_l)\frac{\partial  ue}{\partial y}(x_k,y_l,t_{s+1})\bigg ]\\ \nonumber
 &- \bigg[(\sigma _{i,l}(x_k,y_l)\frac{\partial ^2 v}{\partial x^2}(x_k,y_l,t_{s+1})+ 
\sigma _{i,l,x}(x_k,y_l)\frac{\partial  v}{\partial x}(x_k,y_l,t_{s+1})) + ((\sigma _{i,t}(x_k,y_l)
\frac{\partial ^2 v}{\partial y^2} \\ & (x_k,y_l,t_{s+1})
 +\sigma _{i,t,y}(x_k,y_l)\frac{\partial  v}{\partial y}(x_k,y_l,t_{s+1})\bigg ]= I_1-I_2,\\
\label{discmwb2d}
& \frac{\partial w}{\partial t}(x_k,y_l,t_{s+1})=g(v(x_k,y_l,t_{s+1}),w(x_k.y_l,t_{s+1}).
\end{align}
Using \eqref{discwb2d} at the grid points in \eqref{discmwb2d} and linearize the non-linear terms by treating it explicitly, we get
the following
\begin{align*}
&\sum_{i=1}^{2M}\gamma_{r,s} h_r(x_k) h_s(y_l) = g(v(x_k,y_l,t_{s}),w(x_k,y_l,t_{s})).
\end{align*}

Matrix system of the above equation is given by,
\begin{align}
\label{matwb2d}
H_r H_s \gamma=c,
\end{align}
where $H_r, H_s$ are the Haar matrices and $c^t=(c_{k,l})$, which is given by,
\begin{align}
c_{k.l}=g(v(x_k,y_l,t_{s}),w(x_k,y_l,t_{s})).
\end{align}
Now, at each time step we will calculate the wavelet coefficient $\gamma$ and then from \eqref{wb2d} at the collocation points we will 
calculate the solution $w$. So, now we will use this $w$ to calculate the solution $v$ and $u_e$.

Again, Calculate equations \eqref{ueb2y}, \eqref{ueb2x}, \eqref{ueb1x}, \eqref{ueb1y} and \eqref{vb1t} at 
the collocation points and substitute in \eqref{discmvb2d} and treat non-linear terms explicitly in $v$, we get the following
\begin{align}
\nonumber 
& \sum_{i,j=1}^{2M}\alpha_{i,j} P_{2,i}(x_k) P_{2,j}(y_l) +\bigg (\frac{\partial v}{\partial t}\bigg )(0,y_l,t)- 
\bigg (\frac{\partial v}{\partial t}\bigg )(0,0,t) +\bigg (\frac{\partial v}{\partial t}\bigg )(x_k,0,t)\\ \nonumber
&-\sigma _{e,l}(x_k, y_l) \bigg [ \sum_{m,n=1}^{2M}\beta_{m,n} h_m(x_k) P_{2,n}(y_l) 
 +\frac{\partial ^2 u_e}{\partial x^2}(x_k,0,t) \bigg]-\sigma _{e,l,x}(x_k,y_l) \\ \nonumber
 &\bigg [\sum_{m,n=1}^{2M} \beta_{m,n} P_{1,m}(x_k) P_{2,n}(y_l)+\frac{\partial ue}{\partial x}(x_k,0,t)
 \bigg ]-\sigma _{e,t}(x_k, y_l) \bigg [ \sum_{m,n=1}^{2M}\beta_{m,n} h_m(x_k)\\ \nonumber
& P_{2,n}(y_l)+\frac{\partial ^2 u_e}{\partial x^2}(x_k,0,t) \bigg]-\sigma _{e,t,y_l}(x_k,y_l)
\bigg [\sum_{m,n=1}^{2M}\beta_{m,n} P_{1,m}(x_k) P_{2,n}(y_l) + \frac{\partial ue}{\partial y}\\
&(0,y_l,t)\bigg ]+f(v(x_k,y_l,t_{s}),w(x_k,y_l,t_{s+1})) = I_2
\end{align}

\begin{align}
\nonumber 
&-(\sigma _{i,l}+\sigma _{e,l})(x_k, y_l) \bigg [ \sum_{m,n=1}^{2M}\beta_{m,n} h_m(x_k) P_{2,n}(y_l) 
 +\frac{\partial ^2 u_e}{\partial x^2}(x_k,0,t) \bigg]-(\sigma _{i,l,x}+\sigma _{e,l,x})\\
 &(x_k,y_l)\bigg [\sum_{m,n=1}^{2M} \beta_{m,n} P_{1,m}(x_k) P_{2,n}(y_l)+\frac{\partial ue}{\partial x}(x_k,0,t)
 \bigg ]-(\sigma _{i,t}+\sigma _{e,t})(x_k, y_l) \bigg [ \sum_{m,n=1}^{2M}\\
 &\beta_{m,n} h_m(x_k) P_{2,n}(y_l)+\frac{\partial ^2 u_e}{\partial x^2}(x_k,0,t) \bigg]-\sigma _{e,t,y_l}(x_k,y_l)
\bigg [\sum_{m,n=1}^{2M}\beta_{m,n} P_{1,m}(x_k) P_{2,n}(y_l) + \frac{\partial ue}{\partial y}\\
&(0,y_l,t)\bigg ]+f(v(x_k,y_l,t_{s}),w(x_k,y_l,t_{s+1})) = I_2
\end{align}

The above equation in matrix form at time $t_{s+1}$ can be written as follows:
\begin{align}
\label{matv2d}
K \begin{bmatrix} \alpha & \beta \end{bmatrix} =b,
\end{align}
where $ K=(k_{ij})$ ia matrix of size $8M^2 \times 8M^2$ and $b^t=(b_i)$ is a column vector of size $8M^2 \times 1$.

Now from the above equation we will calculate the wavelet coefficients $\alpha, \beta$ and the obtain the solution $v$ with the use of calculated $w$, at the desired time step.

\section{Analysis of Convergence}
In this section we discuss the convergence analysis for the proposed numerical scheme.

\begin{lemma}  
	If $v(x,y)$ and $w(x,y)$ are Lipschitz continuous on domain $[0,1]^2$, then the wavelet coefficients
	$a_{i_1,i_2}$, $b_{l_1,l_2}$ corresponding to $v$ and $w$ satisfy the inequality 
	\begin{align}
	a_{i_1,i_2} \leq \frac{C}{4m^3},\\
	b_{l_1,l_2} \leq \frac{C}{4m^3}.
	\end{align}
\end{lemma}

Introducing the norm 
\begin{align}
\label{norm}
\parallel u \parallel_X^2 = \parallel v \parallel_2^2 + \parallel u_e \parallel_2^2 + \sum_{k=1}^{d}\parallel w^k \parallel_2^2,
\end{align}
where $ \parallel . \parallel_2 $ is the standard $L^2$-norm. 

Let  $ u =\begin{bmatrix} v & u_e & w \end{bmatrix}$ be the exact solution of the problem and
$ u_H =\begin{bmatrix} v_H & u_{e,h} & w_H \end{bmatrix}$ be the solution approximated by the Haar wavelets. 
The error $(E = u-u_H)$ is given as follows:

\begin{theorem}
	Let $u_H$ be the solution approximated by the Haar wavelet, then
	\begin{align}
	\parallel E \parallel_X ^2 = \parallel u(x,y,t_{n+1})-u_H(x,y,t_{n+1}) \parallel_X ^2 \leq \frac{C^2}{16} 
	\bigg[ \frac{K_1 K_2}{m^2} + \frac{K_3 K_4}{m^2}+ \frac{d}{m^6} \bigg]
	\end{align}
\end{theorem} 
Proof. 
\begin{align*}
\parallel E \parallel_X ^2 &=  \parallel dt \sum_{i_1,i_2=2M+1}^{\infty} a_{i_1,i_2} P_{2,i_1}(x) P_{2,i_2}(y) \parallel_2 ^2+  \parallel \sum_{j_1,j_2=2M+1}^{\infty} b_{j_1,j_2} P_{2,j_1}(x) P_{2,j_2}(y) \parallel_2 ^2, \\ 
& +  \sum_{w=1}^d \parallel dt \sum_{k_1,k_2=2M+1}^{\infty} c_{k_1,k_2}^w H_{k_1}(x) H_{k_2}(y) \parallel_2 ^2 \\
& = dt^2 \int_{0}^{1} \int_{0}^{1} |\sum_{i_1,i_2=2M+1}^{\infty} a_{i_1,i_2} P_{2,i_1}(x) P_{2,i_2}(y) \sum_{r_1,s_1=2M+1}^{\infty} a_{r_1,s_1} P_{2,r_1}(x) P_{2,s_1}(y)| \\
&+  \parallel \sum_{j_1,j_2=2M+1}^{\infty} b_{j_1,j_2} P_{2,j_1}(x) P_{2,j_2}(y) \parallel_2 ^2  P_{2,i_2}(y) \sum_{r_2,s_2=2M+1}^{\infty} b_{r_2,s_2} P_{2,r_2}(x) P_{2,s_2}(y)| \\
&+ \sum_{w=1}^d \parallel dt^2 \int_{0}^{1} \int_{0}^{1} | \sum_{k_1,k_2=2M+1}^{\infty} c_{k_1,k_2}^w H_{k_1}(x) H_{k_2}(y) \sum_{p,q=2M+1}^{\infty} c_{p,q}^k H_{p}(x) H_{q}(y)| \\
& \leq dt^2 \sum_{i_1,i_2=2M+1}^{\infty} |a_{i_1,i_2} |  \sum_{r_1,s_1=2M+1}^{\infty} |a_{r_1,s_1} | \bigg(\int_{0}^{1}|P_{2,i_1}(x) P_{2,r_1}(x)|dx \bigg)  \bigg(\int_{0}^{1}|P_{2,i_2}(y) P_{2,s_1}(y)|dy\bigg)\\
&+ \sum_{j_1,j_2=2M+1}^{\infty} |b_{j_1,j_2} |  \sum_{r_2,s_2=2M+1}^{\infty} |b_{r_2,s_2} | \bigg(\int_{0}^{1}|P_{2,j_1}(x) P_{2,r_2}(x)|dx \bigg)  \bigg(\int_{0}^{1}|P_{2,j_2}(y) P_{2,s_2}(y)|dy\bigg)\\
&+ dt^2 \sum_{w=1}^d \bigg[\sum_{k_1,k_2=2M+1}^{\infty} |c_{k_1,k_2}^k| \sum_{p,q=2M+1}^{\infty} |c_{p,q}^k| \bigg( \int_{0}^{1} |H_{k}(x) H_{p}(x)|dx \bigg) \bigg(\int_{0}^{1} |H_{m}(y) H_{q}(y)|dy\bigg) \bigg]\\
&= dt^2 \sum_{i_1, i_2=2M+1}^{\infty} |a_{i_1, i_2}| \sum_{r_1,r_2=2M+1}^{\infty} |a_{r_1,r_2}| K_{i_1, r_1}K_{i_2, r_2} + \sum_{j_1, j_2=2M+1}^{\infty} |a_{j_1, j_2}| \sum_{s_1, s_2=2M+1}^{\infty} |a_{s_1,s_2}| K_{j_1, s_1}K_{j_2, s_2}  \\
& + dt^2 \sum_{w=1}^d \bigg[\sum_{k_1,k_2=2M+1}^{\infty} |c_{k_1,k_2}^k| \sum_{p,q=2M+1}^{\infty} |c_{p,q}^w| L_{k_1,p} L_{k_2,q} \bigg],\\
&= I_1+I_2 +I_3,
\end{align*}
where
\begin{align*}
&K_{i,r}=\int_{0}^{1}|P_{2,i}(x) P_{2,r}(x)|dx,\\
&L_{j,p}= \int_{0}^{1} |H_{j}(x) H_{p}(x)|dx.
\end{align*}
Let $K_{i_1} = \sup_{r_1} K_{i_1,r_1}$ and $K_{i_2} = \sup_{r_2} K_{i_2, r_2}$, first term $(I_1)$ becomes as follows
\begin{align*}
I_1 & \leq  dt^2 \sum_{i_1, i_2=2M+1}^{\infty} |a_{i_1, i_2}| K_{i_1}K_{i_2} \sum_{r_1,r_2=2M+1}^{\infty} |a_{r_1,r_2}|\\
& \leq  \frac{Cdt^2}{4} \sum_{i_1, i_2=2M+1}^{\infty} |a_{i_1, i_2}| K_{i_1}K_{i_2} \sum_{j=J+1}^ \infty \sum_{r_1,r_2=0}^{2^{j}-1} \frac{1}{m^3} \\
& \leq  \frac{Cdt^2}{4}  \sum_{i_1, i_2=2M+1}^{\infty} |a_{i_1, i_2}| K_{i_1}K_{i_2} \sum_{j=J+1}^ \infty \frac{1}{m} \\
& \leq  \frac{Cdt^2}{4} \sum_{i_1, i_2=2M+1}^{\infty} |a_{i_1, i_2}| K_{i_1}K_{i_2} \frac{1}{m}.
\end{align*}

Now, let $K_1 = \sup_i K_{i_1}$ and $K_2 = \sup_n K_{i_2}$, we get
\begin{align*}
I_1 & \leq  \frac{Cdt^2}{4}  K_1 K_2  \sum_{i_1, i_2=2M+1}^{\infty} |a_{i_1, i_2}| \frac{1}{m} \\
& \leq  \frac{C^2dt^2}{16}  K_1 K_2 \sum_{j=J+1}^ \infty \sum_{i_1, i_2 = 0}^{2^{j}-1} \frac{1}{m^4} \\
I_1& \leq  \frac{C^2dt^2 K_1 K_2}{16m^2}
\end{align*}
Similarly, the second term becomes as 
\begin{align*}
I_2 & \leq  \frac{C^2 K_3 K_4}{16m^2}
\end{align*}
Again, we know by the definition of Haar wavelets ,  
\begin{align*}
\int_{0}^{1} H_{j}(x) H_{p}(x)dx =  2^{-j} =\frac{1}{m}
\end{align*}
\begin{align*}
I_2 &= (dt)^2 \sum_{w=1}^d \bigg[\sum_{k_1,k_2=2M+1}^{\infty} |b_{k_1, k_2}^w|^2 \frac{1}{m^2} \bigg]\\
& \leq \frac{C^2dt^2}{16}\sum_{w=1}^d \bigg(\frac{1}{m^3}\bigg)^2\\
& \leq \frac{dC^2dt^2}{16m^6}
\end{align*}
Therefore,
\begin{align*}
\parallel E \parallel_X ^2 & \leq \frac{C^2dt^2 K_1 K_2}{16m^2} + \frac{C^2 K_3 K_4}{16m^2} + \frac{dC^2dt^2}{16m^6} \\
& \leq \frac{C^2}{16} \bigg[ \frac{K_1 K_2}{m^2} + \frac{K_3 K_4}{m^2}+ \frac{d}{m^6} \bigg]
\end{align*}

\section{Haar wavelet method for the coupled degenerate reaction diffusion PDEs and the ODEs in three dimension}
\begin{align}
\label{mv3d}
&C_m \frac{\partial v}{\partial t}+div(D_e(x)\nabla u_e) + f(v,w)=I_2 \hspace{1.5cm} x \in \Omega , t \in (0,T) \\  
\label{mue3d}
&- div((D_i+)D_e)(x)\nabla u_e)-div(D_i(x)\nabla v)=I_1-I_2; \hspace{0.6cm} x \in \Omega, t \in (0,T)\\  
\label{mw3d}
&\frac{\partial w}{\partial t}-g(v,w)=0 \hspace{1cm} x \in  \Omega,  t \in (0,T)\\  \nonumber
&v(x,0)= v_0(x,0), \hspace{5mm} w(x,0)=w_0(x,0) \hspace{1.8cm}  x \in \Omega\\  \nonumber
&n^T D_{i,e}(x) \nabla v =0 \hspace{ 5.2cm}  x \in \partial \Omega, t \in (0,T). 
\end{align}

Let us write $\frac{\partial ^7 v}{\partial t \partial x^2 \partial y^2 \partial z^2}\big (x,y,z,t\big )$, 
$\frac{\partial ^6 u_e}{\partial x^2 \partial y^2 \partial z^2}\big (x,y,z,t\big )$ and $\frac{\partial w}{\partial t}$ in 
terms of the Haar wavelet as follows:
	
\begin{align}
\label{approxv3d}
\frac{\partial ^7 v}{\partial t \partial x^2 \partial y^2 \partial z^2} (x,y,z,t)&=\sum_{i_1,i_2,i_3=1}^{2M}\alpha_{i_1,i_2,i_3} h_{i_1}(x) h_{i_2}(y) h_{i_3}(z), \hspace{1cm}t\in [t_s, t_{s+1})\\
\label{approxue3d}
\frac{\partial ^6 u_e}{\partial x^2 \partial y^2 \partial z^2} (x,y,z,t)&=\sum_{j_1,j_2,j_3=1}^{2M}\beta_{j_1,j_2,j_3} h_{j_1}(x) h_{j_2}(y) h_{j_3}(z), \hspace{1cm}t\in [t_s, t_{s+1})\\
\label{approxw3d}
\frac{\partial w}{\partial t}(x,y,z,t) &= \sum_{l,m,n=1}^{2*M}\gamma_{l,m,n} h_l(x) h_m(y) h_n(z) , \hspace{0.5cm}t\in [t_s, t_{s+1}).
\end{align}
Integrating equation \eqref{approxv3d} w.r.t  $t$ from $t_s$ to $t$, we will get 
\begin{align}
\label{3dv6}
\frac{\partial ^6 v}{ \partial x^2 \partial y^2  \partial z^2}(x,y,z,t) &= (t-t_s)\sum_{i,j,k=1}^{2M}\alpha_{i,j,k} h_i(x) h_j(y)h_k(z)
+\frac{\partial ^6 v}{ \partial x^2 \partial y^2  \partial z^2} (x,y,z,t_s), \hspace{0.5cm}t\in [t_s, t_{s+1}).
\end{align}

Now, Integrate equation \eqref{3dv6} twice w.r.t  $x$ from $0$ to $x$ also using the boundary conditions, we will obtain the following
\begin{align}
\label{3dv4}
\nonumber
\frac{\partial ^4 v}{ \partial y^2 \partial z^2}(x,y,z,t) &= (t-t_s)\sum_{i,j,k=1}^{2M}\alpha_{i,j,k} P_{2,i}(x) h_j(y)h_k(z)
+\frac{\partial ^4 v}{ \partial y^2  \partial z^2} (x,y,z,t_s) -\frac{\partial ^4 v}{ \partial y^2  \partial z^2} (0,y,z,t_s) \\ &+ \frac{\partial ^4 v}{ \partial y^2  \partial z^2} (0,y,z,t), \hspace{0.5cm}t\in [t_s, t_{s+1}).
\end{align}

Now, Integrate equation \eqref{3dv4} twice w.r.t  $y$ from $0$ to $y$ also using the boundary conditions, we get
\begin{align}
\label{3dvz2}
\nonumber
\frac{\partial ^2 v}{\partial z^2}(x,y,z,t) &= (t-t_s)\sum_{i,j,k=1}^{2M}\alpha_{i,j,k} P_{2,i}(x) P_{2,j}(y)h_k(z) + 
\frac{\partial ^2 v}{\partial z^2}(x,y,z,t_s)-\frac{\partial ^2 v}{\partial z^2}(x,0,z,t_s)-\frac{\partial ^2 v}{\partial z^2}
(0,y,z,t_s)\\
&+ \frac{\partial ^2 v}{\partial z^2}(0,0,z,t_s)+ \frac{\partial ^2 v}{\partial z^2}(0,y,z,t)- \frac{\partial ^2 v}{\partial z^2}
(0,0,z,t)+ \frac{\partial ^2 v}{\partial z^2}(x,0,z,t),
\hspace{0.5cm}t\in [t_s, t_{s+1}).
\end{align}

Similarly, Integrate equation \eqref{3dv4} twice w.r.t  $z$ from $0$ to $z$ also using the boundary conditions, we get	
\begin{align}
\label{3dvy2}
\nonumber
\frac{\partial ^2 v}{\partial y^2}(x,y,z,t) &= (t-t_s)\sum_{i,j,k=1}^{2M}\alpha_{i,j,k} P_{2,i}(x)h_j(y) P_{2,k}(z) + 
\frac{\partial ^2 v}{\partial y^2}(x,y,z,t_s)-\frac{\partial ^2 v}{\partial y^2}(x,y,0,t_s)-\frac{\partial ^2 v}{\partial y^2}
(0,y,z,t_s)\\
&+ \frac{\partial ^2 v}{\partial y^2}(0,y,0,t_s)+ \frac{\partial ^2 v}{\partial y^2}(0,y,z,t)- \frac{\partial ^2 v}{\partial y^2}
(0,y,0,t)+ \frac{\partial ^2 v}{\partial y^2}(x,y,0,t),
\hspace{0.5cm}t\in [t_s, t_{s+1}).
\end{align}

Again, Integrating \eqref{3dv6} twice w.r.t  $z$ from $0$ to $z$ and then twice w.r.t  $y$ from $0$ to $y$ also using the boundary conditions, we get
\begin{align}
\label{3dvx2}
\nonumber
\frac{\partial ^2 v}{\partial x^2}(x,y,z,t) &= (t-t_s)\sum_{i,j,k=1}^{2M}\alpha_{i,j,k} h_i(x) P_{2,j}(y) P_{2,k}(z) + 
\frac{\partial ^2 v}{\partial x^2}(x,y,z,t_s)-\frac{\partial ^2 v}{\partial x^2}(x,0,z,t_s)-\frac{\partial ^2 v}{\partial x^2}(x,y,0,t_s)\\
&+ \frac{\partial ^2 v}{\partial x^2}(x,0,0,t_s)+ \frac{\partial ^2 v}{\partial x^2}(x,0,z,t)- 
\frac{\partial ^2 v}{\partial x^2}(x,0,0,t)+ \frac{\partial ^2 v}{\partial x^2}(x,y,0,t),
\hspace{0.5cm}t\in [t_s, t_{s+1}).
\end{align}

Now, Integrate equation \eqref{approxv3d} twice w.r.t. x , y and then z also using boundary conditions, we will obtain
\begin{align}
\label{3dvt}
\nonumber
\frac{\partial v}{\partial t}(x,y,z,t) &= \sum_{i,j,k=1}^{2M}\alpha_{i,j,k} P_{2,i}(x) P_{2,j}(y) P_{2,k}(z) + \frac{\partial v}{\partial t}(0,y,z,t)- \frac{\partial v}{\partial t}(0,y,0,t)-\frac{\partial v}{\partial t}(0,0,z,t) + \frac{\partial v}{\partial t}(0,0,0,t)\\
& + \frac{\partial v}{\partial t}(x,0,z,t)- \frac{\partial v}{\partial t}(x,0,0,t) + \frac{\partial v}{\partial t}(x,y,0,t), \hspace{0.5cm}t\in [t_s, t_{s+1}).
\end{align}

Now, Integrating the above equation \eqref{3dvt} w.r.t. $t$ from $t_s$ to $t$, we will get 
\begin{align}
\label{v3d}
\nonumber
v(x,y,z,t) &= (t-t_s)  \sum_{i,j,k=1}^{2M}\alpha_{i,j,k} P_{2,i}(x) P_{2,j}(y) P_{2,k}(z) + v(x,y,z,t_s)+ v(0,y,z,t) -v(0,y,z,t_s) - v(0,y,0,t)\\ \nonumber
&+v(0,y,0,t_s)- v(0,0,z,t) + v(0,0,z,t_s) + v(0,0,0,t)- v(0,0,0,t_s) + v(x,0,z,t) - v(x,0,z,t_s)\\
& - v(x,0,0,t) + v(x,0,0,t_s) + v(x,y,0,t) -  v(x,y,0,t_s).
\end{align}

Now, Integrate equation \eqref{approxue3d} twice w.r.t  $x$ from $0$ to $x$ also using the Neumann boundary condition on $u_e$, we will obtain the following
\begin{align*}
\nonumber
\frac{\partial ^4 u_e}{ \partial y^2 \partial z^2}(x,y,z,t) &= (t-t_s)\sum_{j_1,j_2,j_3=1}^{2M}\alpha_{j_1,j_2,j_3} P_{2,j_1}(x) h_{j_2}(y)h_{j_3}(z)
+\frac{\partial ^4 u_e}{ \partial y^2  \partial z^2} (0,y,z,t_s)  \hspace{0.5cm}t\in [t_s, t_{s+1}).
\end{align*}
Integrate the above equation twice w.r.t  $x$ from $0$ to $x$ also using the Neumann boundary condition on $u_e$, we will get
\begin{align}
\label{3duez2}
\frac{\partial ^2 u_e}{\partial z^2}(x,y,z,t) =\sum_{i,j,k=1}^{2M}\alpha_{i,j,k} P_{2,i}(x) P_{2,j}(y)h_k(z) + \frac{\partial ^2 u_e}{\partial z^2}(x,0,z,t)+\frac{\partial ^2 u_e}{\partial z^2}(0,y,z,t)- \frac{\partial ^2 v}{\partial z^2}(0,0,z,t),
\hspace{0.3cm}t\in [t_s, t_{s+1}).
\end{align}

Similarly, first Integrate the equation \eqref{approxue3d} twice w.r.t  $z$ from $0$ to $z$ and then integrate the obtain equation twice w.r.t $x$ also using the boundary conditions, we get	
\begin{align}
\label{3duey2}
\frac{\partial ^2 u_e}{\partial y^2}(x,y,z,t) = \sum_{i,j,k=1}^{2M}\alpha_{i,j,k} P_{2,i}(x)h_j(y) P_{2,k}(z) + \frac{\partial ^2 v}{\partial y^2}(x,y,0,t)+\frac{\partial ^2 v}{\partial y^2}(0,y,z,t)- \frac{\partial ^2 v}{\partial y^2}(0,y,0,t),
\hspace{0.3cm}t\in [t_s, t_{s+1}).
\end{align}
Similarly, first Integrate the equation \eqref{approxue3d} twice w.r.t  $z$ from $0$ to $z$ and then integrate the obtain equation twice w.r.t $y$ from $0$ to $y$ also using the boundary conditions, we get	
\begin{align}
\label{3duex2}
\frac{\partial ^2 u_e}{\partial x^2}(x,y,z,t) &= \sum_{i,j,k=1}^{2M}\alpha_{i,j,k} h_i(x) P_{2,j}(y) P_{2,k}(z) + \frac{\partial ^2 u_e}{\partial x^2}(x,0,z,t)+\frac{\partial ^2 u_e}{\partial x^2}(x,y,0,t)+ \frac{\partial ^2 u_e}{\partial x^2}(x,0,0,t),
\hspace{0.3cm}t\in [t_s, t_{s+1}).
\end{align}
Now, integrate the above equation twice w.r.t. $x$ from $0$ to $x$ also using the boundary conditions, we get	
\begin{align}
u_e(x,y,z)=  \sum_{i,j,k=1}^{2M}\beta_{i,j,k} h_i(x) P_{2,j}(y) P_{2,k}(z) + u_e(x,0,z,t)- u_e(0,0,z,t) + u_e(x,y,0,t) -u_e(0,y,0,t) -u_e(x,0,0,t) +u_e(0,0,0,t). 
\end{align}

Again, Integrate \eqref{approxw3d} w.r.t $t$ from $t_s$ to $t$, we acquire
\begin{align}
\label{w3d}
w(x,y,z,t) = (t-t_s)\sum_{l,m,n=1}^{2M}\gamma_{l,m,n} h_{l}(x) h_{m}(y) h_{n}(z)+ w(x,y,z,t_s)
\end{align}

To find the solution at the collocation points, we have to discretized the equation \eqref{mv3d} - \eqref{mw3d} when 
$t \rightarrow t_{s+1}$.
The discrete form is as follows:
\begin{align}
\nonumber
\label{discmv3d}
& \frac{\partial v}{\partial t}(x_{k_1},y_{k_2},z_{k_3},t_{s+1})+ \bigg[\sigma _{l,e}(x_{k_1},y_{k_2},z_{k_3})\frac{\partial ^2 u_e}{\partial x^2}(x_{k_1},y_{k_2},z_{k_3},t_{s+1})+
\sigma _{l,e,x}(x_{k_1},y_{k_2},z_{k_3})\frac{\partial  u_e}{\partial x}(x_{k_1},y_{k_2},z_{k_3},t_{s+1}) \\ \nonumber
&+\sigma _{t,e}(x_{k_1},y_{k_2},z_{k_3})\frac{\partial ^2 u_e}{\partial y^2}(x_{k_1},y_{k_2},z_{k_3},t_{s+1})
+\sigma _{t,e,y}(x_{k_1},y_{k_2},z_{k_3})\frac{\partial  u_e}{\partial y}(x_{k_1},y_{k_2},z_{k_3},t_{s+1}) +\sigma _{t,e}(x_{k_1},y_{k_2},z_{k_3})\\ \nonumber
&\frac{\partial ^2 u_e}{\partial z^2}(x_{k_1},y_{k_2},z_{k_3},t_{s+1})
+\sigma _{t,e,z}(x_{k_1},y_{k_2},z_{k_3})\frac{\partial  u_e}{\partial z}(x_{k_1},y_{k_2},z_{k_3},t_{s+1})\bigg ]
+f(v(x_{k_1},y_{k_2},z_{k_3},t_{s+1}),w(x_{k_1},y_{k_2},z_{k_3},\\ &t_{s+1})) = I_2
\end{align}

\begin{align}
\nonumber
& -\bigg[(\sigma _{l,i}+\sigma _{l,e})(x_{k_1},y_{k_2},z_{k_3})\frac{\partial ^2 u_e}{\partial x^2}(x_{k_1},y_{k_2},z_{k_3},t_{s+1})+
(\sigma _{l,i,x}+\sigma _{l,e,x})(x_{k_1},y_{k_2},z_{k_3})\frac{\partial  u_e}{\partial x}(x_{k_1},y_{k_2},z_{k_3},t_{s+1}) \\ \nonumber
&+((\sigma _{t,i}+\sigma _{t,e})(x_{k_1},y_{k_2},z_{k_3})\frac{\partial ^2 u_e}{\partial y^2}(x_{k_1},y_{k_2},z_{k_3},t_{s+1})
+(\sigma _{t,i,y}+\sigma _{t,e,y})(x_{k_1},y_{k_2},z_{k_3})\frac{\partial  u_e}{\partial y}(x_{k_1},y_{k_2},z_{k_3},t_{s+1}) +(\sigma _{t,i}+\sigma _{t,e})(x_{k_1},y_{k_2},z_{k_3})\\ \nonumber
&\frac{\partial ^2 u_e}{\partial z^2}(x_{k_1},y_{k_2},z_{k_3},t_{s+1})
+(\sigma _{t,i,z}+\sigma _{t,e,z})(x_{k_1},y_{k_2},z_{k_3})\frac{\partial  u_e}{\partial z}(x_{k_1},y_{k_2},z_{k_3},t_{s+1})\bigg ] -  \bigg[\sigma _{l,i}(x_{k_1},y_{k_2},z_{k_3})\frac{\partial ^2 v}{\partial x^2}(x_{k_1},y_{k_2},z_{k_3},t_{s+1})+ \\ \nonumber
&\sigma _{l,i,x}(x_{k_1},y_{k_2},z_{k_3})\frac{\partial v}{\partial x}(x_{k_1},y_{k_2},z_{k_3},t_{s+1}) \\ \nonumber
&+\sigma _{t,i}(x_{k_1},y_{k_2},z_{k_3})\frac{\partial ^2 v}{\partial y^2}(x_{k_1},y_{k_2},z_{k_3},t_{s+1})
+\sigma _{t,e,y}(x_{k_1},y_{k_2},z_{k_3})\frac{\partial  v}{\partial y}(x_{k_1},y_{k_2},z_{k_3},t_{s+1}) +\sigma _{t,i}(x_{k_1},y_{k_2},z_{k_3})\\ 
\label{discmue3d}
&\frac{\partial ^2 v}{\partial z^2}(x_{k_1},y_{k_2},z_{k_3},t_{s+1})
+\sigma _{t,i,z}(x_{k_1},y_{k_2},z_{k_3})\frac{\partial v}{\partial z}(x_{k_1},y_{k_2},z_{k_3},t_{s+1})\bigg ] = I_1-I_2,
\end{align}

\begin{align}
\label{discmw3d}
& \frac{\partial w}{\partial t} (x_{k_1},y_{k_2},z_{k_3},t_{s+1}) = g(v(x_{k_1},y_{k_2},z_{k_3},t_{s+1}),w(x_{k_1},y_{k_2},z_{k_3},t_{s+1}).
\end{align}
Using \eqref{approxw3d} at the grid points in \eqref{discmw3d} and linearize the non-linear terms by treating it explicitly, we obtain the following
\begin{align*}
&\sum_{l,m,n=1}^{2M}\gamma_{l,m,n} h_l(x_{k_1}) h_m(y_{k_2}) h_n(z_{k_3}) = g(v(x_{k_1},y_{k_2},z_{k_3},t_{s}),w(x_{k_1},y_{k_2},z_{k_3},t_{s})).
\end{align*}
Matrix system of the above equation is given by,
\begin{align}
\label{matw2d}
H_l H_m H_n \gamma=c,
\end{align}
where $H_l, H_m, H_n$ are the Haar matrices and $c^t=(c_{k_1, k_2, k_3})$, which is given by,
\begin{align}
c_{{k_1, k_2, k_3}}= g(v(x_{k_1},y_{k_2},z_{k_3},t_{s}),w(x_{k_1},y_{k_2},z_{k_3},t_{s})).
\end{align}	

Now, at each time step we will calculate the wavelet coefficient $\beta$ and then from \eqref{w3d} at the collocation points we
will calculate the solution $w$. So, now we will use this $w$ to calculate the solutions $u_e$ and $v$.

Again, Calculate equations \eqref{ueb2y}, \eqref{ueb2x}, \eqref{ueb1x}, \eqref{ueb1y} and \eqref{vb1t} at 
the collocation points and substitute in \eqref{discmvb2d} and treat non-linear terms explicitly in $v$, we get the following

Again, Calculate equations \eqref{3dvy2}, \eqref{3dvx2}, \eqref{3dvx1}, \eqref{3dvy1}, \eqref{3dvt},  \eqref{3duey2}, \eqref{3duex2}, \eqref{3duex1}, \eqref{3duey1}, at the collocation points and substitute in equations \eqref{discmv3d} and \eqref{discmue3d} and treat non-linear terms explicitly in $v$, we will get the following matrix system at time $t_{s+1}$ :
\begin{align}
\label{matv2d}
K \begin{bmatrix} \alpha & \beta \end{bmatrix} =b,
\end{align}
where $ K=(k_{ij})$ ia matrix of size $8M^3 \times 8M^3$ and $b^t=(b_i)$ is a column vector of size $8M^3 \times 1$.

Now from the above equation we will calculate the wavelet coefficients $\alpha, \beta$ and obtain the solutions $u_e$ and $v$ with the
use of calculated $w$, at the desired time step.

\section{Numerical Result and Discussions}
We solve all the examples using above developed haar wavelet method and calculate the absolute error also. Grid validation test
or resolution level test has been done for all the problems and here we are presenting for some of the problems. From grid 
validation we observe that resolution level $J=4$ in two dimension is good enough to calculate the solution. We use the GMRES solver to solve the linear system of equations.

\begin{figure}[h]
	\centering
	\begin{subfigure}[t]{0.4\textwidth}
		\centering
		\includegraphics[width=1.5\textwidth]{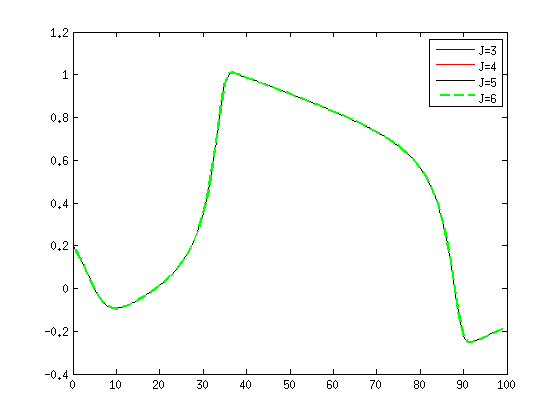}
		\caption{}
		\label{gridValid_1D_v}
	\end{subfigure}\hfill
	\begin{subfigure}[t]{0.4\textwidth}
		\centering
		\includegraphics[width=1.5\textwidth]{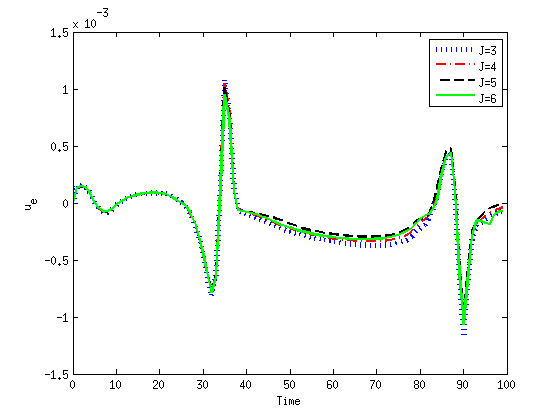}
		\caption{}
		\label{gridValid_1D_ue}
	\end{subfigure}\hfill
	\caption{Grid validation test for, (a)$v$, (b) $u_e$}
	\label{gridValid_1D}
\end{figure}

\textbf{Example $1$.}  We consider the one dimensional degenerate coupled PDEs and the ODE having homogeneous Neumann boundary as follows:
\begin{align*}
&C_m \frac{\partial v}{\partial t}+\frac{d}{dx}(D_e(x)\frac{d}{dx} u_e) +  v(v-0.1)(1-v)-w=I^e_{app}, \hspace{0.6cm}0\leq x \leq 1,  0\leq t \leq T\\ 
&- \frac{d}{dx}(D_i+D_e)(x)\frac{d}{dx} u_e)-\frac{d}{dx}(D_i(x)\frac{d}{dx} v)=I^i_{app}-I^e_{app}, \hspace{0.6cm}0\leq x \leq 1,  0\leq t \leq T\\
&\frac{\partial w}{\partial t}=v-2w \hspace{1cm} 0\leq x \leq 1,  0\leq t \leq T\\
&v(x,0)=0.2, \hspace{5mm} w(x,0)=0.2 \hspace{1.8cm}  0\leq x \leq 1,\\
&D_{i,e}(x)\frac{d u_{i,e}}{d x}(0,t)=0, \hspace{0.5cm} D_{i,e}(x)\frac{du_{i,e}}{d x}(1,t)=0 \hspace{1cm}  0\leq t \leq T. 
\end{align*}
where, 
Resolution level test for the proposed Haar wavelet method has been presented in Fig. \ref{gridValid_1D}.

\begin{figure}[h]
	\centering
	\begin{subfigure}[t]{0.4\textwidth}
		\centering
		\includegraphics[width=1.3\textwidth]{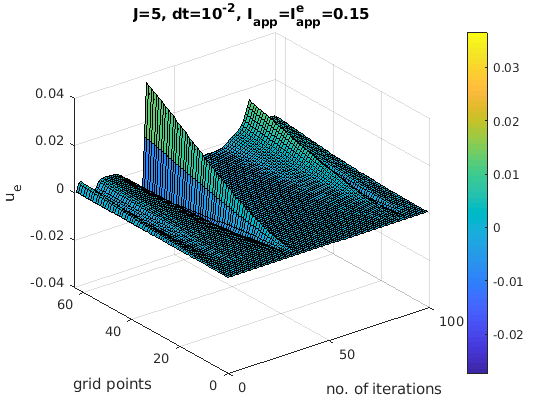}
		\caption{{Haar wavelet solution for $u_e$}}
		\label{surf1Due}
	\end{subfigure}\hfill
	\begin{subfigure}[t]{0.4\textwidth}
		\includegraphics[width=1.3\textwidth]{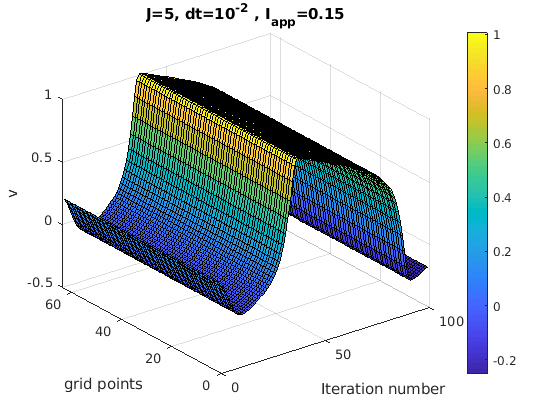}
		\caption{{Haar wavelet solution for $v$}}
		\label{surf1Dv}
	\end{subfigure}\hfill
	\caption{Haar wavelet solution}
	\label{surf1D_fhn}
\end{figure}

\begin{figure}
\centering
	\includegraphics[width=0.6\textwidth]{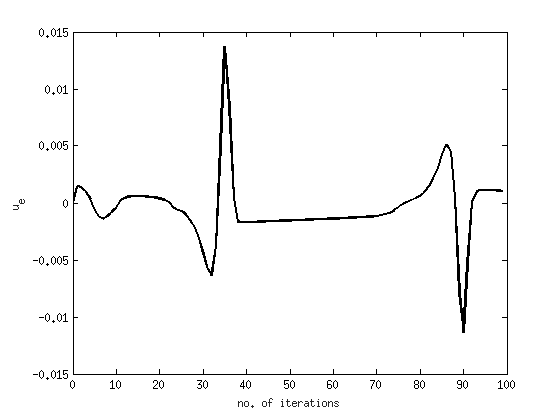}
	\caption{\textbf{Haar wavelet solution for $u_e$ at point (0.5, 0.5)}}
	\label{ue1d}
\end{figure}

Pointwise absolute error for different time steps is shown in Table \ref{erro1d_v} and \ref{erro1d_ue}.
Solution at $dt=10^{-5}, J=7$ is taken as the reference solution. From the Table \ref{erro1d_v} and \ref{erro1d_ue}.
it can be seen clearly that absolute error for $v$ and $u_e$ decreases significantly with the smaller time step size.
The Haar wavelet solution for $v$ and $u_e$ at the grid points are given in Fig. \ref{surf1D_fhn} and \ref{ue1d}. 

\begin{table}
	\begin{center}
		\begin{tabular}{ |c c  c  c | } 
			\hline
			x & & absolute error & \\
			\hline 
		    & $dt=10^{-2}$	& $dt=10^{-3}$ & $dt =10^{-4}$\\
			\hline
			0.0234 & $2.83 \times10^{-2}$  & $2.0 \times10^{-3}$ & $1.34 \times 10^{-4}$\\ 
			
			0.1172 & $2.83 \times10^{-2}$ & $2.0 \times10^{-3}$ & $1.34 \times 10^{-4}$\\ 
			
			0.2266 & $2.82 \times10^{-2}$ & $2.0 \times10^{-3}$ & $1.33 \times 10^{-4}$\\ 
			
			0.3828 & $2.81 \times10^{-2}$ & $2.0 \times10^{-3}$ & $1.32 \times 10^{-4}$\\ 
			
			0.5391 & $2.81 \times10^{-2}$ & $2.0 \times10^{-3}$ & $1.32 \times 10^{-4}$\\
				
	    	0.7734 & $2.81 \times10^{-2}$ & $2.0 \times10^{-3}$ & $1.31 \times 10^{-4}$\\
	    	
	    	0.9297 & $2.81 \times10^{-2}$ & $2.0 \times10^{-3}$ & $1.3 \times 10^{-4}$\\	
			\hline
		\end{tabular}
		\caption{Absolute Error for $v$ at different points of the domain when $J=5$ and T=0.5}
		\label{erro1d_v}
	\end{center}
\end{table} 

\begin{table}
	\begin{center}
		\begin{tabular}{ |c c  c  c | } 
			\hline
			x & & absolute error & \\
			\hline 
			& $dt=10^{-2}$	& $dt=10^{-3}$ & $dt =10^{-4}$\\
			\hline
			0.0234 & $2.2048 \times10^{-5}$  & $2.952 \times10^{-6}$ & $1.952 \times 10^{-6}$\\ 
			
			0.1172 & $1.8499 \times10^{-4}$ & $1.011 \times10^{-6}$ & $2.011 \times 10^{-6}$\\ 
			
			0.2266 & $5.1482 \times10^{-4}$ & $3.822 \times10^{-6}$ & $2.1780 \times 10^{-6}$\\ 
			
			0.3828 & $1.1 \times10^{-3}$ & $1.4401 \times10^{-5}$ & $1.599 \times 10^{-6}$\\ 
			
			0.5391 & $1.7 \times10^{-3}$ & $2.5823 \times10^{-5}$ & $1.77 \times 10^{-7}$\\
			
			0.7734 & $2.4 \times10^{-3}$ & $4.8658 \times10^{-3}$ & $6.58 \times 10^{-7}$\\
			
			0.9297 & $2.8 \times10^{-3}$ & $6.3827 \times10^{-5}$ & $8.27 \times 10^{-7}$\\	
			\hline
		\end{tabular}
		\caption{Absolute Error for $u_e$ at different points of the domain when $J=5$ and T=0.5}
		\label{erro1d_ue}
	\end{center}
\end{table}

\paragraph{\textbf{Example $2$.}} We consider the two dimensional degenerate coupled PDEs and the ODE having homogeneous
Neumann boundary as follows:
\begin{align*}
&C_m \frac{\partial v}{\partial t}+ \nabla. (D_e(x,y)\nabla u_e) + v(v-0.1)(1-v)-kw =I^e_{app}, \hspace{0.5cm}  0\leq x,y \leq 1,  0\leq t \leq T\\ 
&-\nabla ((D_i+D_e)(x,y)\nabla u_e)-\nabla (D_i(x,y)\nabla v)=I^i_{app}-I^e_{app} \hspace{0.6cm} 0\leq x,y \leq 1,  0\leq t \leq T\\
&\frac{\partial w}{\partial t}-g(v,w) =0 \hspace{1cm} 0\leq x,y \leq 1,  0\leq t \leq T\\
\text{with initial condition}\\
&v(x,y,0)= v_0(x,y,0), \hspace{5mm} w(x,y,0) =w_0(x,0) \hspace{1.8cm} 0\leq x,y \leq 1\\
\text{and Neuman boundary conditions}\\
&D(x,y)\frac{d u_{i,e}}{d x}(0,y,t)=0,\hspace{2cm} 0\leq t \leq T
&D(x,y)\frac{d u_{i,e}}{d x}(1,y,t) =0 \hspace{2cm}, 0\leq t \leq T\\
&D(x,y)\frac{d u_{i,e}}{d y}(x,0,t)=0,\hspace{2cm} 0\leq t \leq T\\
&D(x,y)\frac{d u_{i,e}}{d y}(x,1,t) =0, \hspace{2cm}, 0\leq t \leq T
\end{align*}

\begin{figure}[h]
	\centering
	\begin{subfigure}[t]{0.4\textwidth}
		\centering
		\includegraphics[width=1.4\textwidth]{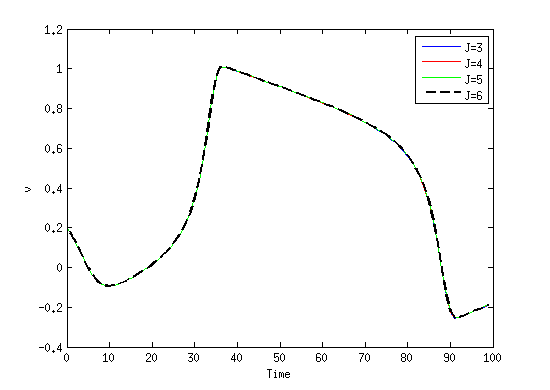}
		\caption{}
		\label{gridValid_2D_v}
	\end{subfigure}\hfill
	\begin{subfigure}[t]{0.4\textwidth}
		\centering
		\includegraphics[width=1.3\textwidth]{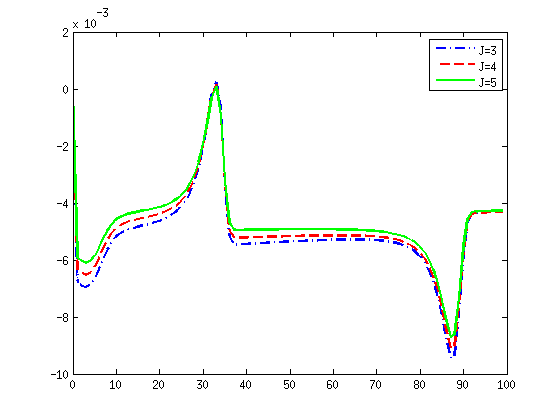}
		\caption{}
		\label{gridValid_2D_ue}
	\end{subfigure}\hfill
	\caption{Grid validation test for, (a)$v$, (b) $u_e$}
	\label{gridValid_2D}
\end{figure}

First of all the grid validation of the proposed algorithm for this problem is presented in Fig. \ref{gridValid_2D} which clearly 
shows the accuracy of the solution at the different resolution level. So, resolution level $J=4$ is good enough to calculate the 
results.

Error for different time steps is shown in Table \ref{error_2d}. Solution at $dt=10^{-5}, J=5$ is taken as the reference solution. 
From the Table \ref{error_2d}, it can be seen clearly that error decreases seriously with the smaller time step size.

\begin{figure}[h]
	\centering
	\begin{subfigure}[t]{0.4\textwidth}
		\centering
		\includegraphics[width=1.5\textwidth]{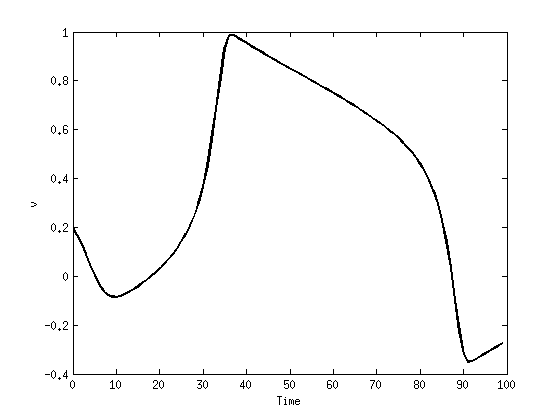}
		\caption{{Haar wavelet solution for $v$}}
		\label{v_fhn_2d}
	\end{subfigure}\hfill
	\begin{subfigure}[t]{0.4\textwidth}
		\centering
		\includegraphics[width=1.5\textwidth]{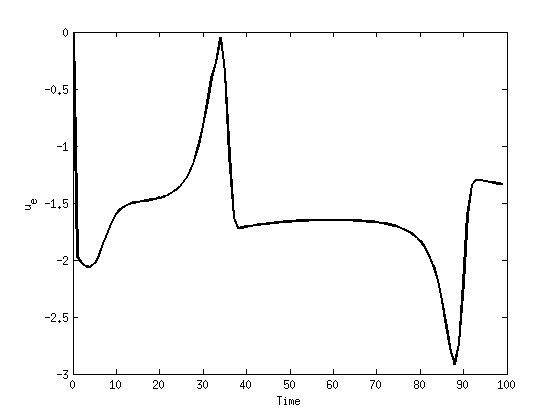}
		\caption{{Haar wavelet solution for $u_e$}}
		\label{ue_fhn2d}
	\end{subfigure}\hfill
	\caption{Haar wavelet solution in two dimension}
	\label{surf2D_fhn}
\end{figure}

\begin{table}[h]
	\begin{center}
		\begin{tabular}{ |c c  c  c | } 
			\hline 
			& $dt=10^{-2}$	& $dt=10^{-3}$ & $dt =10^{-4}$\\
			\hline
			$L^{\infty}(ue)$ error & $1.17 \times10^{-4}$  & $1.5 \times10^{-5}$ & $1.0 \times 10^{-6}$\\ 
				\hline
			$L^{\infty}(v)$ error & $1.17 \times10^{-4}$  & $1.5 \times10^{-5}$ & $1.0 \times 10^{-6}$\\ 
				\hline	
		\end{tabular}
		\caption{Error for $u_e$ and $v$ at T=1 with reference solution when $J=5, dt=10^{-5}$.}
		\label{error_2d}
	\end{center}
\end{table}

\begin{figure}[h]
	\centering
	\begin{subfigure}[t]{0.4\textwidth}
		\centering
		\includegraphics[width=1.5\textwidth]{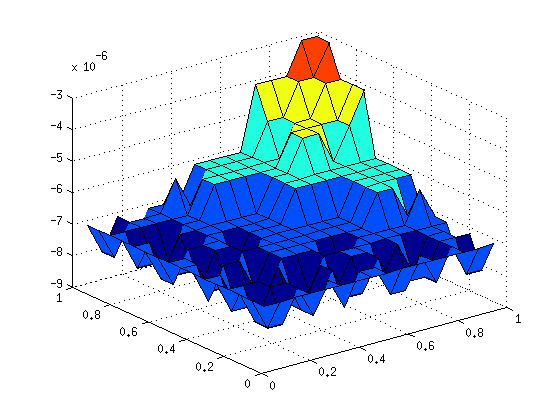}
		\caption{{Haar wavelet solution error for $v$}}
		\label{v_fhn_2d_error}
	\end{subfigure}\hfill
	\begin{subfigure}[t]{0.4\textwidth}
		\centering
		\includegraphics[width=1.5\textwidth]{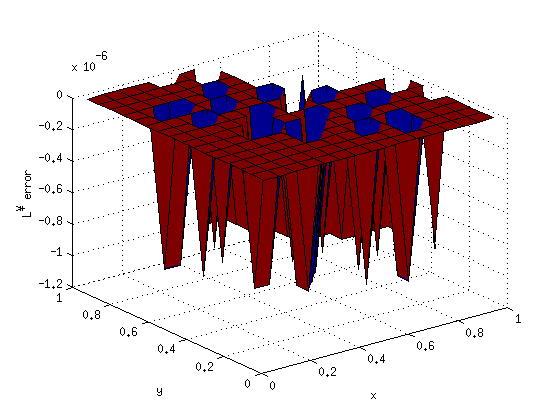}
		\caption{{Haar wavelet solution error for $u_e$}}
		\label{ue_fhn2d_error}
	\end{subfigure}\hfill
	\caption{Haar wavelet solution error in 2 dimension}
	\label{surf2D_fhn_error}
\end{figure}

\paragraph{Example $3$}
\begin{align*}
	C_m \frac{\partial v}{\partial t} + \nabla .(D_e(x,y,z)\nabla u_e) + v(v-0.1)(1-v)-kw &=I^e_{app} \hspace{1.5cm}  0\leq x,y \leq 1,
	0\leq t \leq T\\
 - \nabla ((D_i+D_e)(x,y,z)\nabla u_e)-\nabla (D_i(x,y,z)\nabla v)&=I^i_{app}-I^e_{app} \hspace{0.6cm} 0\leq x,y \leq 1,  0\leq t \leq T\\
	\frac{\partial w}{\partial t}&=v-2w \hspace{2cm} 0\leq x,y,z \leq 1,  0\leq t \leq T\\
	v(x,y,z,0)= 0.2, \hspace{5mm} w(x,y,z,0)&=0.2 \hspace{1cm} 0\leq x,y,z \leq 1	
	\end{align*}
with Neumann boundary conditions on $v$ and $u_e$. 
	
where $d_{11}= 1.2 \times 10^{-3}, d_{22}=2.5562\times 10^{-4}, d_{33}=2.5562 \times 10^{-4}$.

\begin{figure}[h]
	\centering
	\begin{subfigure}[t]{0.4\textwidth}
		\centering
		\includegraphics[width=1.2\textwidth]{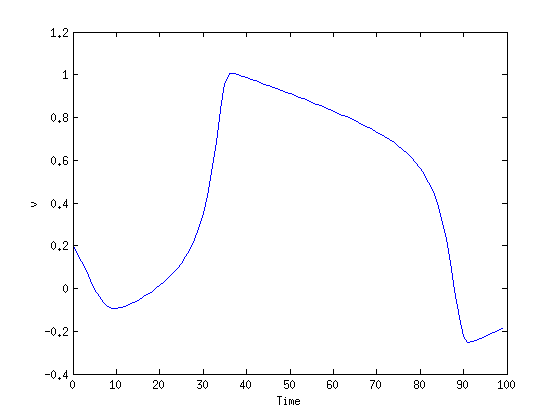}
		\caption{{Haar wavelet solution for $v$}}
		\label{v_fhn_3d}
	\end{subfigure}\hfill
	\begin{subfigure}[t]{0.4\textwidth}
		\centering
		\includegraphics[width=1.2\textwidth]{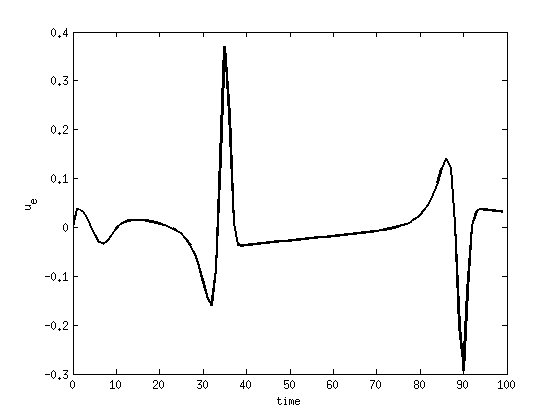}
		\caption{{Haar wavelet solution for $u_e$}}
		\label{ue_fhn3d}
	\end{subfigure}\hfill
	\caption{Haar wavelet solution in 3 dimension}
	\label{plot3d_fhn}
\end{figure}

The solution of this system for $v$ and $u_e$ using Haar wavelet method is presented in Fig. \ref{plot3d_fhn}.

\begin{figure}
\centering
	\includegraphics[width=0.8\textwidth]{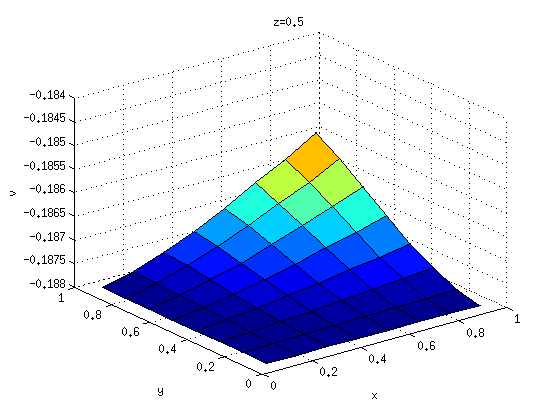}
		\caption{{Haar wavelet solution for $v$}}
		\label{v_fhn_3d_surf}
	\end{figure}
	
Remark:	The examples presented above are vastly applicable to the field of cardiac electrophysiology. 

\section{Conclusion}
A Haar Wavelet Method for a class of for the coupled degenerate reaction diffusion PDEs and the ODEs having non-linear source with Neumann boundary has been proposed. The method is both simple and easy to implement in two and three dimensions. 
Convergence analysis has also been done to ensure the stability and accuracy. Model problems have been successfully solved. Numerical error reduces with the increase in time step size or resolution level. Problems with clinical relevance have also been successfully dealt with.

\end{document}